\documentclass[11pt,a4paper]{article}
\usepackage[T1]{fontenc}
\usepackage[utf8]{inputenc}
\usepackage{lmodern}
\usepackage[margin=1.05in,includefoot]{geometry}
\usepackage{amsmath,amssymb,amsthm,mathtools}
\usepackage{enumitem}
\usepackage{tikz}
\usepackage[hidelinks,unicode]{hyperref}
\input{glyphtounicode}
\hypersetup{pdftitle={A Proof of the Balister--Gy\H{o}ri--Schelp Conjecture on Perfect Matchings with Prescribed Differences},
  pdfauthor={Aryeh Lev Zabokritskiy (Yohananov)},pdflang={en}}

\theoremstyle{plain}
\newtheorem{theorem}{Theorem}[section]
\newtheorem{lemma}[theorem]{Lemma}
\newtheorem{proposition}[theorem]{Proposition}
\newtheorem{corollary}[theorem]{Corollary}
\theoremstyle{definition}
\newtheorem{definition}[theorem]{Definition}
\newtheorem{example}[theorem]{Example}
\newtheorem{remark}[theorem]{Remark}
\numberwithin{equation}{section}

\DeclareMathOperator{\Pf}{Pf}
\DeclareMathOperator{\tr}{tr}
\DeclareMathOperator{\sgn}{sgn}
\newcommand{\F}{\mathbb{F}_2}
\newcommand{\Q}{\mathbb{Q}}
\newcommand{\Z}{\mathbb{Z}}
\newcommand{\nut}{\nu_2}
\newcommand{\mm}{\mathfrak{m}}
\newcommand{\form}[2]{\langle #1,#2\rangle}
\newcommand{\cG}{\mathcal{G}}

\title{A Proof of the Balister--Gy\H{o}ri--Schelp Conjecture on Perfect Matchings with Prescribed Differences}
\author{
  \normalsize Aryeh Lev Zabokritskiy (Yohananov)\\[0.5em]
  \small Department of Computer Science, MIGAL--Galilee Research Institute,\\
  \small Kiryat Shmona, Israel, and Department of Computer Science, Tel-Hai\\
  \small University of Kiryat Shmona and the Galilee, Kiryat Shmona, Israel\\[0.25em]
  \small\href{mailto:yuhanalev@telhai.ac.il}{\texttt{yuhanalev@telhai.ac.il}}\\[0.5em]
  \small ORCID: \href{https://orcid.org/0000-0003-3151-6192}{0000-0003-3151-6192}
}
\date{}

\begin{document}
\maketitle

\begin{abstract}
We prove the conjecture of Balister, Gy\H{o}ri and Schelp on perfect matchings
with prescribed differences. Let $n\ge2$ and let $d_1,\dots,d_{2^{n-1}}$
be nonzero vectors in $\F^n$, with repetitions allowed. The space $\F^n$
can be partitioned into pairs $\{x_i,y_i\}$ with $x_i+y_i=d_i$ for every
$i$ if and only if $\sum_i d_i=0$.
\end{abstract}

\section{Introduction}\label{sec:intro}

Let $n\ge 2$, let $V\coloneqq\F^n$, and put $M\coloneqq 2^{n-1}$, so that
$|V|=2M$.  Given
nonzero vectors $d_1,\dots,d_M\in V$, not necessarily distinct, we ask whether
$V$ can be partitioned into $M$ pairs in such a way that the two elements of
the $i$-th pair differ by $d_i$.  Formally, a \emph{realization} of the list
$(d_1,\dots,d_M)$ is a sequence $(E_1,\dots,E_M)$ of pairwise disjoint subsets
of $V$ such that $E_i=\{x_i,x_i+d_i\}$ for some $x_i\in V$.  The sets
$E_1,\dots,E_M$ then form a perfect matching of $V$ whose edges have the
prescribed differences $d_1,\dots,d_M$, and each index $i$ has its own edge
even when several of the vectors $d_i$ coincide.

A realization can exist only if $d_1+\dots+d_M=0$.  Indeed, $d_i$ is the sum
of the two elements of $E_i$, so $d_1+\dots+d_M$ is the sum of all elements of
$V$, and this sum is zero because every coordinate equals one on exactly
$2^{n-1}$ vectors and $2^{n-1}$ is even.  Balister, Gy\H{o}ri and Schelp
conjectured that this condition is also sufficient~\cite{BGS}.
Previous results cover several special cases~\cite{GK,Kovacs,FourEight,TwoHole}.
We prove the conjecture in full.

\begin{theorem}\label{thm:main}
Let $n\ge 2$ and $M=2^{n-1}$.  A list $(d_1,\dots,d_M)$ of nonzero vectors of
$\F^n$ has a realization if and only if $d_1+\dots+d_M=0$.
\end{theorem}

By the connection described by Yohananov and Essayag~\cite{YE} (see also~\cite{HKRS}),
our theorem settles a conjecture of Zhang, Etzion and Yaakobi~\cite[Conjecture~25]{ZEY}:
for any list of $2^{n-1}$ linear forms on the $n$ message bits, their values can
be recovered from pairwise disjoint sets of coordinates of the binary simplex code,
one set per form.

As in~\cite{TwoHole}, the proof rests on a signed count.
For a list with sum zero, fix a total order on $V$ and choose an ordered
basis in which each coordinate is the last nonzero coordinate of an even
number of the $d_i$, counted with multiplicity. Such a basis always exists.
These choices give each realization a sign. We prove that for every such
choice the sum of the signs is an odd multiple of $2^{M-1}$. In particular,
it is nonzero, so at least one realization exists.

In Section~\ref{sec:certificate} we express the signed count as a Pfaffian
coefficient.  Section~\ref{sec:traces} gives a trace formula for the
coefficients of the logarithm of the normalized Pfaffian.  In
Section~\ref{sec:trees} we group the terms by binary trees, and in
Section~\ref{sec:parity} we complete the proof by showing that an odd number
of trees contribute to the coefficient for the full list.
Section~\ref{sec:conclusion} asks whether a realizing matching can be found
efficiently.

\section{The signed count of realizations and its valuation}
\label{sec:certificate}

Throughout the paper $n\ge 2$, $V=\F^n$ and $M=2^{n-1}$, as in the
Introduction; we put $N\coloneqq 2^n=|V|$, and $d=(d_1,\dots,d_M)$ is a list
of nonzero vectors of $V$.  The elements of
$[M]\coloneqq\{1,\dots,M\}$ are called \emph{labels}; thus $d_i$ is the prescribed
difference with label $i$, and a realization of $d$ is a sequence
$(E_1,\dots,E_M)$ of pairwise disjoint sets $E_i=\{x_i,x_i+d_i\}$, one for each
label.  We write $V^*$ for the space of linear functionals
$V\to\F$, and $(-1)^c$ denotes $1$ for $c=0$ and $-1$ for $c=1$ in $\F$.

\subsection{Realizations as terms of a Pfaffian}\label{subsec:pfaffian}

The Pfaffian provides a sign for each perfect matching; for the identities
below see, for example,~\cite[Chapter~7]{Godsil}.  Fix a total order on $V$.  A \emph{perfect
matching} of $V$ is a partition of $V$ into sets of size two.  Let
$A=(A_{x,y})_{x,y\in V}$ be a skew-symmetric matrix with zero diagonal and
entries in a commutative ring.  Its Pfaffian is
\begin{equation}\label{eq:pfaffian}
  \Pf(A)=\sum_{\mu}\sgn(\mu)\prod_{\substack{\{x,y\}\in\mu\\ x<y}}A_{x,y},
\end{equation}
where $\mu$ runs over the perfect matchings of $V$, and $\sgn(\mu)\in\{1,-1\}$
is the sign of the permutation that lists the elements of $V$ pair by pair,
the smaller element of each pair first, compared with the listing of $V$ in
increasing order.  The sign $\sgn(\mu)$ depends on $\mu$
and on the order of $V$ only, because listing the pairs in another order
changes the permutation by an even permutation.  Moreover,
\begin{equation}\label{eq:pf-facts}
  \Pf(A)^2=\det A
  \qquad\text{and}\qquad
  \Pf(Q^{\mathsf T}AQ)=\det(Q)\,\Pf(A)
\end{equation}
for every square matrix $Q$ over the same ring.

To see a realization as one term of~\eqref{eq:pfaffian} we use skew-symmetric
matrices whose nonzero entries lie at the positions $(x,x+v)$ for a fixed
vector $v$.  For a nonzero $v\in V$ and a functional $\eta\in V^*$ with
$\eta(v)=1$, let $J(v,\eta)$ be the matrix with rows and columns indexed by
$V$ and entries
\[
  J(v,\eta)_{x,y}\coloneqq
  \begin{cases}
    (-1)^{\eta(x)} & \text{if } y=x+v,\\
    0 & \text{otherwise.}
  \end{cases}
\]
It is skew-symmetric, because
$J(v,\eta)_{x+v,x}=(-1)^{\eta(x)+\eta(v)}=-J(v,\eta)_{x,x+v}$, and its diagonal
is zero, because $v\ne 0$.

Now let $\eta_1,\dots,\eta_M\in V^*$ be functionals with $\eta_i(d_i)=1$, put
$J_i\coloneqq J(d_i,\eta_i)$, and let $t_1,\dots,t_M$ be independent
variables, one for each label.  The \emph{signed count of realizations} of $d$
with respect to $\eta_1,\dots,\eta_M$ is the integer
\begin{equation}\label{eq:Z}
  Z\coloneqq[t_1t_2\cdots t_M]\,\Pf(t_1J_1+\dots+t_MJ_M),
\end{equation}
where $[t_1t_2\cdots t_M]$ denotes the coefficient of the monomial
$t_1t_2\cdots t_M$.  A change of the order of $V$ amounts to replacing the
matrix $t_1J_1+\dots+t_MJ_M$ by $Q^{\mathsf T}(t_1J_1+\dots+t_MJ_M)Q$ for a
permutation matrix $Q$, so it multiplies $Z$ by $\det Q=\pm 1$, by the second
identity in~\eqref{eq:pf-facts}. Each realization contributes one sign to $Z$.

\begin{lemma}[Coefficient certificate]\label{lem:certificate}
We have
\[
  Z=\sum_{(E_1,\dots,E_M)}\varepsilon(E_1,\dots,E_M),\qquad
  \varepsilon(E_1,\dots,E_M)\coloneqq\sgn(\mu)\prod_{i=1}^{M}(-1)^{\eta_i(x_i)},
\]
where the sum runs over all realizations of $d$. For each realization
$(E_1,\dots,E_M)$, let $\mu$ be its underlying perfect matching and
let $x_i$ be the smaller element of $E_i$. We call
$\varepsilon(E_1,\dots,E_M)\in\{1,-1\}$ the \emph{sign} of the
realization.  In particular, if $Z\ne 0$, then $d$ has a realization.
\end{lemma}

\begin{proof}
By~\eqref{eq:pfaffian},
\[
  \Pf(t_1J_1+\dots+t_MJ_M)
  =\sum_{\mu}\sgn(\mu)\prod_{\substack{\{x,y\}\in\mu\\ x<y}}
   \Bigl(\sum_{i=1}^{M}t_i\,(J_i)_{x,y}\Bigr).
\]
Expanding the product over the $M$ pairs of $\mu$ amounts to choosing a label
for each pair.  The resulting monomial is $t_1t_2\cdots t_M$ exactly when the
choice is a bijection between the pairs of $\mu$ and $[M]$.  Let $E_i$ be the
pair that receives the label $i$, and write $E_i=\{x_i,y_i\}$ with $x_i<y_i$.
The corresponding term is $\sgn(\mu)\prod_{i=1}^{M}(J_i)_{x_i,y_i}$.  It is
nonzero if and only if $y_i=x_i+d_i$ for every $i$, that is, if and only if
$(E_1,\dots,E_M)$ is a realization of $d$, and in that case it equals
$\varepsilon(E_1,\dots,E_M)$, because $(J_i)_{x_i,x_i+d_i}=(-1)^{\eta_i(x_i)}$.
Each realization arises in this way exactly once, from the matching
$\{E_1,\dots,E_M\}$ and the choice of the label $i$ for the pair $E_i$.
\end{proof}

Equal prescribed differences still carry distinct labels in this count.
The certificate is sufficient for the existence of a realization, but whether
$Z\ne 0$ depends on the functionals $\eta_i$: realizations with opposite signs
cancel in $Z$.  When all the $d_i$ lie in the nonzero coset $a+H$ of a linear
hyperplane $H$, a matching pairs each $x\in H$ with $a+\pi(x)$ for a
permutation $\pi$ of $H$.  Thus this is exactly Hall's
prescribed-differences problem on $H$~\cite{Hall}: the displacements
$\pi(x)-x$ must be the $d_i+a$, counted with multiplicity.  The \emph{BSolution}
algorithm of Yohananov and Yaakobi constructs a realization in this binary
case~\cite[Algorithm~7 and Theorem~41]{YY}.  An alternative proof
in~\cite[Section~3]{TwoHole} uses the exact $2$-adic valuation of a
coefficient of a group circulant determinant in place of the Pfaffian.
Minors of that circulant treat the case with exactly two differences in the
hyperplane.

\begin{remark}[A necessary condition on the functionals]\label{rem:balance}
If $\eta_1+\dots+\eta_M\ne 0$, then $Z=0$.  To see this, let $a\in V$ and let
$Q$ be the permutation matrix of the translation $x\mapsto x+a$.  Then
$(Q^{\mathsf T}J(v,\eta)Q)_{x,y}=J(v,\eta)_{x+a,y+a}=(-1)^{\eta(a)}J(v,\eta)_{x,y}$,
and $\det Q=1$, because for $a\ne 0$ the translation is a product of $M$
transpositions and $M$ is even.  By~\eqref{eq:pf-facts}, the Pfaffian
in~\eqref{eq:Z} is unchanged when each $t_i$ is replaced by
$(-1)^{\eta_i(a)}t_i$.  Comparing the coefficients of $t_1t_2\cdots t_M$ gives
$Z=(-1)^{(\eta_1+\dots+\eta_M)(a)}Z$ for every $a\in V$, and the claim follows.
\end{remark}

\subsection{Even flags and the choice of the functionals}\label{subsec:flag}

We now choose the functionals $\eta_i$.  By Remark~\ref{rem:balance} the
signed count can be nonzero only if their sum vanishes.  We obtain such
functionals from an even flag for the list $d$, defined below; the same flag
is used again in Section~\ref{sec:parity}, where the parity argument pairs
labels of equal height.

A \emph{flag} of $V$ is a chain of subspaces
$\{0\}=H_0\subset H_1\subset\dots\subset H_n=V$ with $\dim H_j=j$.  The set
$H_j\setminus H_{j-1}$ is its $j$-th \emph{layer}, for $1\le j\le n$; the
layers partition $V\setminus\{0\}$.  We call the flag \emph{even} for the list
$d$ if every layer contains $d_i$ for an even number of labels $i$.

\begin{lemma}[Existence of an even flag]\label{lem:flag}
Every list $d=(d_1,\dots,d_M)$ of nonzero vectors of $V$ has an even flag.
\end{lemma}

\begin{proof}
We construct $H_0\subset H_1\subset\dots\subset H_n$ so that each $H_j$
contains $d_i$ for an even number of labels $i$; the layers then have the same
property.  The subspace $H_0=\{0\}$ contains no $d_i$.  Suppose that $H_j$ has
been constructed for some $j<n$.  Since $M$ is even, the number of labels $i$
with $d_i\notin H_j$ is even.  These labels are distributed among the
$2^{n-j}-1$ cosets of $H_j$ different from $H_j$.  The number of these cosets
is odd, so at least one of them, say $x+H_j$, contains $d_i$ for an even
number of labels $i$.  Then $H_{j+1}\coloneqq H_j\cup(x+H_j)$ is a subspace of
dimension $j+1$ with the required property.
\end{proof}

The proof uses two properties of the list: the vectors $d_i$ are nonzero, and
their number $M$ is even.

For the rest of the paper we fix an even flag for $d$ and a basis
$e_1,\dots,e_n$ of $V$ \emph{adapted} to it, that is, $e_j\in
H_j\setminus H_{j-1}$ for every $j$; thus $H_j$ is spanned by $e_1,\dots,e_j$.
Let $e_1^*,\dots,e_n^*\in V^*$ be the dual basis.  For a nonzero $v\in V$ the
\emph{height} $h(v)$ is the index of the layer that contains $v$;
equivalently, $h(v)$ is the largest $j$ with $e_j^*(v)=1$.
We record two consequences of these definitions for later use, in
Section~\ref{subsec:reduction} and in Lemma~\ref{lem:pair}:
\begin{align}
  &e_m^* \text{ vanishes on } H_j \text{ whenever } m>j; \label{eq:vanish}\\
  &\text{if } v,v'\ne 0 \text{ have the same height } j,
   \text{ then } v+v'\in H_{j-1}. \label{eq:sameheight}
\end{align}
The second holds because $H_{j-1}$ has index two in $H_j$.  For the labels we
put $\eta_i\coloneqq e^*_{h(d_i)}$, so that $\eta_i(d_i)=1$, and
$J_i\coloneqq J(d_i,\eta_i)$, and we call $h(d_i)$ the \emph{height of the
label} $i$.  From now on $Z$ denotes the
signed count~\eqref{eq:Z} formed with these functionals.  Since the flag is
even, each $e_j^*$ occurs an even number of times among
$\eta_1,\dots,\eta_M$, and therefore
\begin{equation}\label{eq:balance}
  \eta_1+\dots+\eta_M=0 .
\end{equation}

\begin{example}\label{ex:running}
Let $n=3$ and write the vectors of $V=\F^3$ as binary words $x_1x_2x_3$.  The
list $d=(100,\,010,\,001,\,111)$ has sum zero.  The flag
\[
  \{000\}\subset\{000,110\}\subset\{000,110,100,010\}\subset V
\]
is even for $d$: its three layers contain $d_i$ for $0$, $2$ and $2$ labels.
The basis $e_1=110$, $e_2=100$, $e_3=001$ is adapted to it.  The labels $1,2$
have height $2$ and the labels $3,4$ have height $3$, so
$\eta_1=\eta_2=e_2^*$ and $\eta_3=\eta_4=e_3^*$, where $e_2^*(x)=x_1+x_2$ and
$e_3^*(x)=x_3$.  The list $d$ has exactly eight realizations.  Indeed, a translation of $V$
maps realizations to realizations and permutes the four possible pairs
$E_1=\{x,x+100\}$ transitively, so it suffices to count the realizations with
$E_1=\{000,100\}$.  For these, $E_2$ is $\{001,011\}$ or $\{101,111\}$, and
each choice leaves exactly one possibility for $E_3$ and $E_4$:
\[
\begin{aligned}
  &\bigl(\{000,100\},\{001,011\},\{110,111\},\{010,101\}\bigr),\\
  &\bigl(\{000,100\},\{101,111\},\{010,011\},\{001,110\}\bigr).
\end{aligned}
\]
We return to the signs of the eight realizations after
Theorem~\ref{thm:valuation}.

For comparison, the flag $\{000\}\subset\{000,100\}\subset\{000,100,010,110\}\subset V$
is not even for $d$: its layers contain $d_i$ for $1$, $1$ and $2$ labels.
With the basis $100,010,001$ the functionals attached to the four labels are
the coordinate functionals $x_1$, $x_2$, $x_3$, $x_3$.  Their sum is nonzero, so the corresponding signed
count vanishes by Remark~\ref{rem:balance}: four of the eight realizations
receive the sign $1$ and four the sign $-1$.
\end{example}

\subsection{The valuation theorem}\label{subsec:valuation}

For a nonzero rational number $q$ let $\nut(q)$ denote its $2$-adic valuation,
that is, the integer $r$ for which $q=2^{r}a/b$ with odd integers $a$ and
$b$, and put $\nut(0)\coloneqq+\infty$.  Then $\nut(qq')=\nut(q)+\nut(q')$ and
$\nut(q+q')\ge\min\{\nut(q),\nut(q')\}$, with equality when
$\nut(q)\ne\nut(q')$.  The main assertion of the paper is the following.

\begin{theorem}[Valuation of the signed count]\label{thm:valuation}
Let $n\ge 2$, and let $d=(d_1,\dots,d_M)$ be a list of nonzero vectors of $V$
with $d_1+\dots+d_M=0$.  For every even flag for $d$, every basis adapted to it
and every total order of $V$, the signed count $Z$ of realizations of $d$ with
respect to the functionals $\eta_i=e^*_{h(d_i)}$ satisfies
\[
  \nut(Z)=M-1 .
\]
In particular, $Z\ne 0$.
\end{theorem}

In Example~\ref{ex:running} we have $M=4$, so the theorem gives $\nut(Z)=3$.
Since $Z$ is the sum of the signs of the eight realizations, these signs are
equal and $Z=\pm 8$; Example~\ref{ex:trees} confirms the value $\pm 8$ without
using Theorem~\ref{thm:valuation}.

\begin{proof}[Proof of Theorem~\ref{thm:main}, assuming Theorem~\ref{thm:valuation}]
The necessity of the condition $d_1+\dots+d_M=0$ was shown in the
Introduction.  Conversely, let $d_1+\dots+d_M=0$.  By Lemma~\ref{lem:flag} the
list $d$ has an even flag; choose a basis adapted to it.  Then $Z\ne 0$ by
Theorem~\ref{thm:valuation}, and $d$ has a realization by
Lemma~\ref{lem:certificate}.
\end{proof}

For $n=1$ the zero-sum condition is not necessary: the list consists of the single
vector $d_1=1$, which is realized by $E_1=\F$.

\subsection{Reduction to two propositions}\label{subsec:reduction}

We normalize the Pfaffian to have constant term one and take its logarithm.
The coefficients of the logarithm, indexed by nonempty sets of labels, express
the signed count up to sign as a sum over set partitions.
The contribution of each partition with at least two blocks has $2$-adic
valuation at least $M$.
For the partition with a single block, a parity count of binary trees gives
an odd multiple of $2^{M-1}$ when the prescribed differences sum to zero.

\paragraph{The reference matrix.}
Put $J_0\coloneqq J(e_n,e_n^*)$.  Each row of $J_0$ has exactly one nonzero
entry, in the column $x+e_n$ of the row $x$.  Hence the only perfect matching
that contributes to $\Pf(J_0)$ is $\{\{x,x+e_n\}:x\in V\}$, and
$\Pf(J_0)=\pm 1$.  We define the polynomial
\[
  F\coloneqq\frac{\Pf(J_0+t_1J_1+\dots+t_MJ_M)}{\Pf(J_0)}\in\Z[t_1,\dots,t_M].
\]
Its constant term is $1$.  Its coefficient of $t_1t_2\cdots t_M$ is
$\pm Z$: in the expansion~\eqref{eq:pfaffian} of the numerator, a term of
degree $M$ in the variables takes a variable from each of the $M$ pairs of
$\mu$ and therefore contains no entry of $J_0$.  Consequently
\begin{equation}\label{eq:ZF}
  \nut(Z)=\nut\bigl([t_1t_2\cdots t_M]\,F\bigr).
\end{equation}
The matrix $J_0$ serves two purposes.  It gives $F$ the constant term $1$, so
that the logarithm of $F$ is defined.  Its functional $e_n^*$ vanishes on
$H_{n-1}$ by~\eqref{eq:vanish}, and this property is used in
Lemma~\ref{lem:pair}(ii).

\paragraph{The logarithmic coefficients.}
The coefficient in~\eqref{eq:ZF} is that of a squarefree monomial, so we may
compute in the ring
\[
  R\coloneqq\Q[t_1,\dots,t_M]/(t_1^2,\dots,t_M^2).
\]
For $S\subseteq[M]$ put $t_S\coloneqq\prod_{i\in S}t_i$.  The elements $t_S$
form a basis of $R$ over $\Q$, and $t_St_{S'}$ equals $t_{S\cup S'}$ if $S$
and $S'$ are disjoint and $0$ otherwise.  For $f\in R$ we write $[t_S]f$ for
the coefficient of $t_S$ in $f$; the passage from a polynomial to its image in
$R$ preserves the coefficients of the squarefree monomials.  Let $\mm$ be the
ideal of $R$ spanned by the $t_S$ with $S\ne\emptyset$; it satisfies
$\mm^{M+1}=0$.  For $g\in\mm$ the finite sums
\[
  \log(1+g)\coloneqq\sum_{m=1}^{M}\frac{(-1)^{m-1}}{m}\,g^m,
  \qquad
  \exp(g)\coloneqq\sum_{m=0}^{M}\frac{g^m}{m!}
\]
satisfy the usual identities $\exp(\log(1+g))=1+g$ and
$\log\bigl((1+g)(1+g')\bigr)=\log(1+g)+\log(1+g')$: these are identities of
formal power series, and they may be evaluated at elements of $\mm$ because
$\mm^{M+1}=0$.  The image of $F$ in $R$
lies in $1+\mm$, and we define the \emph{logarithmic coefficients} $L_S\in\Q$
by
\[
  \log F=\sum_{\emptyset\ne S\subseteq[M]}L_S\,t_S .
\]
The coefficient of $t_{[M]}$ in $F$ is recovered from the $L_S$ by the
exponential formula, which takes the following form in $R$.

\begin{lemma}[Expansion over set partitions]\label{lem:partitions}
We have
\[
  [t_{[M]}]\,F=\sum_{\pi}\ \prod_{S\in\pi}L_S ,
\]
where $\pi$ runs over the partitions of the set $[M]$ into nonempty blocks.
\end{lemma}

\begin{proof}
Put $\Lambda\coloneqq\log F$, so that $F=\exp(\Lambda)=\sum_{m}\Lambda^m/m!$.
The coefficient of $t_{[M]}$ in $\Lambda^m$ is the sum of the products
$L_{S_1}\cdots L_{S_m}$ over all sequences $(S_1,\dots,S_m)$ of pairwise
disjoint nonempty sets with union $[M]$.  Every partition of $[M]$ into $m$
blocks arises from exactly $m!$ such sequences, and this factor cancels the
denominator $m!$.
\end{proof}

The two propositions concern the valuations of the logarithmic coefficients.

\begin{proposition}[Proper label sets]\label{prop:proper}
If $S\subseteq[M]$ and $1\le|S|<M$, then $\nut(L_S)\ge|S|$.
\end{proposition}

\begin{proposition}[The full label set]\label{prop:full}
If $d_1+\dots+d_M=0$, then $\nut(L_{[M]})=M-1$.
\end{proposition}

Proposition~\ref{prop:proper} is proved in Section~\ref{sec:trees}.
Proposition~\ref{prop:full} is proved in Sections~\ref{sec:trees}
and~\ref{sec:parity}; it is the only place where the hypothesis
$d_1+\dots+d_M=0$ and the evenness of the flag are used.

\begin{proof}[Proof of Theorem~\ref{thm:valuation}, assuming
Propositions~\ref{prop:proper} and~\ref{prop:full}]
Consider the sum in Lemma~\ref{lem:partitions}.  The partition with the single
block $[M]$ contributes $L_{[M]}$, which has valuation $M-1$ by
Proposition~\ref{prop:full}.  Every other partition $\pi$ has at least two
blocks, each of them a nonempty proper subset of $[M]$, so by
Proposition~\ref{prop:proper} its contribution has valuation at least
$\sum_{S\in\pi}|S|=M$.  Hence $\nut([t_{[M]}]F)=M-1$, and the theorem follows
from~\eqref{eq:ZF}.
\end{proof}

\section{Logarithmic coefficients as sums of traces}\label{sec:traces}

It remains to prove Propositions~\ref{prop:proper} and~\ref{prop:full}.  Both
rest on a formula that expresses $L_S$ through the matrices
$J_0^{-1}J_i$ (Lemma~\ref{lem:trace}).  These matrices belong to a family
of matrices that is closed under products, whose commutation relations are
given by a bilinear form and whose traces lie in $\{0,N,-N\}$
(Lemma~\ref{lem:calculus}); we describe this family first.

For $a\in V$ and $\beta\in V^*$ let $W(a,\beta)$ be the matrix with rows and
columns indexed by $V$ and entries
\[
  W(a,\beta)_{x,y}\coloneqq
  \begin{cases}
    (-1)^{\beta(x)} & \text{if } y=x+a,\\
    0 & \text{otherwise.}
  \end{cases}
\]
Equivalently, $W(a,\beta)$ acts on functions $f\colon V\to\Q$ by
$(W(a,\beta)f)(x)=(-1)^{\beta(x)}f(x+a)$.  A \emph{signed translation} is a
matrix of the form $W(a,\beta)$ or $-W(a,\beta)$, and the pair
$(a,\beta)\in V\times V^*$ is its \emph{class}.  In this notation
$J(v,\eta)=W(v,\eta)$.  On $V\times V^*$ we consider the form
\[
  \form{(a,\beta)}{(c,\gamma)}\coloneqq\beta(c)+\gamma(a)\in\F .
\]
It is bilinear and symmetric, and $\form{\kappa}{\kappa}=0$ for every
$\kappa\in V\times V^*$.

\begin{lemma}[Signed translations]\label{lem:calculus}
\begin{enumerate}[label=\textup{(\roman*)},leftmargin=*]
\item $W(a,\beta)\,W(c,\gamma)=(-1)^{\gamma(a)}\,W(a+c,\beta+\gamma)$.
Consequently a product of signed translations is a signed translation, and its
class is the sum of the classes of the factors.
\item If $X$ and $Y$ are signed translations with classes $\kappa$ and
$\kappa'$, then $XY=(-1)^{\form{\kappa}{\kappa'}}\,YX$.
\item $W(0,0)$ is the identity matrix $I$, and $\tr W(a,\beta)=0$ for
$(a,\beta)\ne(0,0)$.  Hence the trace of a signed translation is $N$, $-N$ or
$0$, and it is $\pm N$ exactly when the class is $(0,0)$.
\end{enumerate}
\end{lemma}

\begin{proof}
(i) For a function $f$ and $x\in V$,
\[
  \bigl(W(a,\beta)W(c,\gamma)f\bigr)(x)
  =(-1)^{\beta(x)}\bigl(W(c,\gamma)f\bigr)(x+a)
  =(-1)^{\beta(x)+\gamma(x+a)}f(x+a+c),
\]
and $\beta(x)+\gamma(x+a)=(\beta+\gamma)(x)+\gamma(a)$.
(ii) By~(i), $W(a,\beta)W(c,\gamma)$ and $W(c,\gamma)W(a,\beta)$ are the
matrix $W(a+c,\beta+\gamma)$ multiplied by $(-1)^{\gamma(a)}$ and
$(-1)^{\beta(c)}$ respectively.
(iii) The diagonal entries of $W(a,\beta)$ vanish for $a\ne 0$, and
$\tr W(0,\beta)=\sum_{x\in V}(-1)^{\beta(x)}$, which is $0$ for $\beta\ne 0$
because a nonzero functional takes each of its two values $N/2$ times.
\end{proof}

By Lemma~\ref{lem:calculus}(i), $J_0^2=(-1)^{e_n^*(e_n)}W(0,0)=-I$.  Hence
$J_0^{-1}=-J_0$, and for every label $i$ the matrix
\[
  B_i\coloneqq J_0^{-1}J_i=-J_0J_i
\]
is a signed translation of class
\begin{equation}\label{eq:classes}
  b_i\coloneqq(e_n+d_i,\ e_n^*+\eta_i)\in V\times V^* .
\end{equation}
For $S\subseteq[M]$ put $b_S\coloneqq\sum_{i\in S}b_i$.  An \emph{ordering} of
a finite set $S$ is a sequence $w=(i_1,\dots,i_k)$ in which every element of
$S$ occurs exactly once.  For an ordering $w=(i_1,\dots,i_k)$ of a set
$S\subseteq[M]$ we write
\[
  B_w\coloneqq B_{i_1}B_{i_2}\cdots B_{i_k};
\]
by Lemma~\ref{lem:calculus}(i) this is a signed translation of class $b_S$.

\begin{lemma}[Trace formula]\label{lem:trace}
Let $S\subseteq[M]$ be a set with $k\ge 1$ elements.  Then
\[
  L_S=\frac{(-1)^{k-1}}{2k}\sum_{w}\tr(B_w),
\]
where $w$ runs over the $k!$ orderings of $S$.
\end{lemma}

\begin{proof}
Let $X\coloneqq t_1B_1+\dots+t_MB_M$.  Since
$J_0+t_1J_1+\dots+t_MJ_M=J_0(I+X)$, the first identity
in~\eqref{eq:pf-facts}, applied to this matrix and to $J_0$, gives
$F^2=\det(I+X)$.  We read this identity in $R$, where the entries of $X$ lie
in $\mm$, so that $X^{M+1}=0$, and take logarithms:
\[
  2\log F=\log\det(I+X)=\tr\log(I+X)
  =\sum_{m=1}^{M}\frac{(-1)^{m-1}}{m}\tr(X^m).
\]
Here the logarithm of the matrix $I+X$ is defined by the same finite series,
\[
  \log(I+X)\coloneqq\sum_{m=1}^{M}\frac{(-1)^{m-1}}{m}\,X^m,
\]
which gives the last equality.  The second equality is the standard relation between determinant
and trace; we verify it at the end of the proof.

Now
$X^m=\sum t_{i_1}\cdots t_{i_m}B_{i_1}\cdots B_{i_m}$, the sum being over all
sequences $(i_1,\dots,i_m)$ of labels.  In $R$ the monomial
$t_{i_1}\cdots t_{i_m}$ vanishes unless the labels are distinct, and it equals
$t_S$ exactly when $(i_1,\dots,i_m)$ is an ordering of $S$, which requires
$m=k$.  Comparing the coefficients of $t_S$ gives the formula.

It remains to verify the second equality.  Let $s$ be a variable.  The
polynomials in $s$ with coefficients in $\mm$ form an ideal of $R[s]$ whose
$(M+1)$-st power is zero, so the same finite series define $\log(1+g)$ for
every element $g$ of this ideal, as well as $\log(I+sX)$.  Termwise
differentiation gives
$\tfrac{d}{ds}\log(1+g)=\tfrac{dg}{ds}\,(1+g)^{-1}$ and
$\tfrac{d}{ds}\tr\log(I+sX)=\tr\bigl((I+sX)^{-1}X\bigr)$, because the terms
neglected in these two formulas are multiples of $\tfrac{dg}{ds}\,g^{M}$ and of
$X^{M+1}$, which vanish.  The element $g=\det(I+sX)-1$ lies in the ideal, and
by Jacobi's formula for the derivative of a determinant, $\log\det(I+sX)$ has
the same derivative $\tr\bigl((I+sX)^{-1}X\bigr)$.  Thus the two sides of
$\log\det(I+sX)=\tr\log(I+sX)$ are polynomials in $s$ with the same derivative
and the same value $0$ at $s=0$.  Since $R$ contains $\Q$, they are equal,
and $s=1$ gives the second equality.
\end{proof}

\section{Summation over orderings by binary trees}\label{sec:trees}

Let $S$ be a nonempty set of $k$ labels. By Lemma~\ref{lem:trace},
the coefficient $L_S$ is, up to the factor
$\pm 1/(2k)$, the sum of the traces of the $k!$ products $B_w$.  Two orderings
that differ by the exchange of two consecutive segments give products
that are equal or opposite, according to Lemma~\ref{lem:calculus}(ii).  We
group the orderings according to binary trees; the orderings in one group
differ by exchanges of this kind.  The sum of the products $B_w$ over the orderings attached to
one tree is then zero or a power of two times a single signed translation
(Lemma~\ref{lem:treesum}).  This gives Proposition~\ref{prop:proper} and
reduces Proposition~\ref{prop:full} to a count of trees modulo two
(Proposition~\ref{prop:fullformula}).

\subsection{Trees, forests and orderings}\label{subsec:trees}

\begin{definition}\label{def:tree}
Let $A$ be a finite set whose size is a power of two.  The \emph{trees on $A$}
and their \emph{nodes}, which are subsets of $A$, are defined recursively.  If
$|A|=1$, there is exactly one tree on $A$, and its only node is $A$.  If
$|A|\ge 2$, a tree on $A$ is an unordered pair $T=\{T',T''\}$ of trees on two
disjoint sets $A'$ and $A''$ with $A'\cup A''=A$ and $|A'|=|A''|$; its nodes
are $A$ and the nodes of $T'$ and of $T''$.  The sets $A'$ and $A''$ are the
\emph{children} of the node $A$.  A node of size two is a \emph{cherry}.
\end{definition}

Thus every node of size at least two has two children, which are nodes of half
its size, and a tree is determined by its set of nodes.  For each element $i$
of $A$ the nodes containing $i$ form a chain with one node of each size
$1,2,4,\dots,|A|$ (the \emph{chain property}); in particular $i$ lies in
exactly one cherry when $|A|\ge 2$.

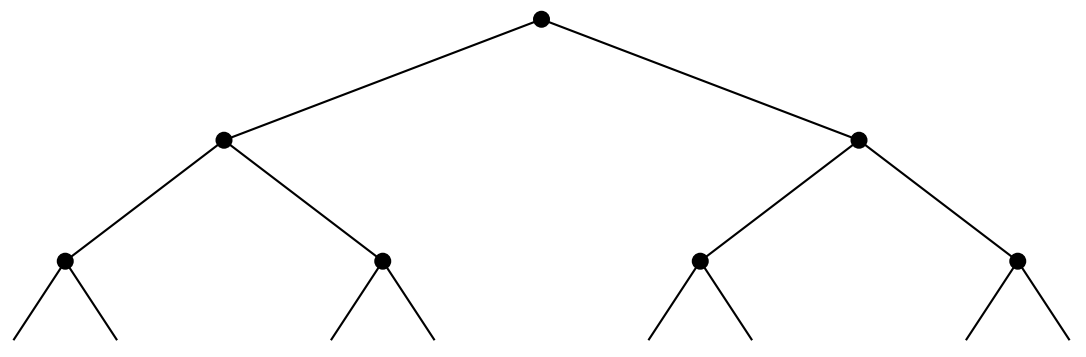
\begin{figure}[ht]
\centering
\begin{tikzpicture}[x=1.05cm,y=0.8cm]
  \foreach \k in {1,...,8} \node[font=\small] (l\k) at (\k,0) {$i_{\k}$};
  \foreach \k/\a/\b in {1/1/2,2/3/4,3/5/6,4/7/8}
    { \coordinate (a\k) at ({(\a+\b)/2},1); \draw (l\a) -- (a\k) -- (l\b); }
  \coordinate (b1) at (2.5,2); \coordinate (b2) at (6.5,2);
  \draw (a1) -- (b1) -- (a2); \draw (a3) -- (b2) -- (a4);
  \coordinate (c) at (4.5,3); \draw (b1) -- (c) -- (b2);
  \foreach \p in {a1,a2,a3,a4,b1,b2,c} \fill (\p) circle (1.6pt);
\end{tikzpicture}
\caption{The tree induced by an ordering $(i_1,\dots,i_8)$.  Each dot
represents the node formed by the labels below it.  The four lowest dots are
the cherries $\{i_1,i_2\}$, $\{i_3,i_4\}$, $\{i_5,i_6\}$ and $\{i_7,i_8\}$.
Exchanging the two halves below any of the seven dots changes the ordering and
leaves the tree unchanged, so $2^7$ orderings induce this tree.}
\label{fig:tree}
\end{figure}

An ordering $w$ of $A$ \emph{induces} a tree $T(w)$ on $A$.  If $|A|=1$, it is
the unique tree on $A$.  If $|A|\ge 2$, split $w$ into its first half $w'$ and
its second half $w''$, and put $T(w)\coloneqq\{T(w'),T(w'')\}$.  Conversely, if
$T=\{T',T''\}$ is a tree on $A$ with $|A|\ge 2$, then the orderings inducing $T$
are exactly the concatenations $w'w''$ and $w''w'$ in which $w'$ induces $T'$
and $w''$ induces $T''$.  Figure~\ref{fig:tree} shows the tree induced by an
ordering of eight labels.

\begin{lemma}\label{lem:count}
Every tree on $A$ is induced by exactly $2^{|A|-1}$ orderings of $A$.
Consequently the number of trees on $A$ is $|A|!/2^{|A|-1}$.
\end{lemma}

\begin{proof}
We use induction on $|A|$; the case $|A|=1$ is clear.  Let $T=\{T',T''\}$ be a
tree on $A$ with $|A|\ge 2$.  By the description of the orderings inducing $T$
given above and by induction, their number is
$2\cdot 2^{|A|/2-1}\cdot 2^{|A|/2-1}=2^{|A|-1}$.  Since each of the $|A|!$ orderings of $A$ induces exactly one tree,
the second statement follows.
\end{proof}

Recall that $b_A=\sum_{i\in A}b_i$ is the class of a product of the matrices
$B_i$ for $i\in A$, in any order. A tree $T$ on $A\subseteq[M]$ is called \emph{good} if every
node of $T$ of size at least two, with children $A'$ and $A''$, satisfies
\[
  \form{b_{A'}}{b_{A''}}=0 .
\]
By Lemma~\ref{lem:calculus}(ii), this means that, for
every such node, a product of the matrices $B_i$ over the labels of one child,
taken in any order, commutes with such a product over the labels of the other
child.  For a tree $T$ on $A$ we write
\[
  \Sigma_T\coloneqq\sum_{w\,:\,T(w)=T}B_w
\]
for the sum of the products $B_w$ over the orderings $w$ of $A$ that induce
$T$.

\begin{lemma}[Sum over the orderings inducing a tree]\label{lem:treesum}
Let $T$ be a tree on a set $A\subseteq[M]$.
\begin{enumerate}[label=\textup{(\roman*)},leftmargin=*]
\item If $T$ is good, then $B_w$ is the same matrix for all orderings $w$
inducing $T$.  We denote it by $B_T$; it is a signed translation of class
$b_A$, and $\Sigma_T=2^{|A|-1}B_T$.
\item If $T$ is not good, then $\Sigma_T=0$.
\end{enumerate}
\end{lemma}

\begin{proof}
We use induction on $|A|$; the case $|A|=1$ is clear.  Let $T=\{T',T''\}$,
where $T'$ and $T''$ are trees on the children $A'$ and $A''$ of $A$.  The
description of the orderings inducing $T$ given before Lemma~\ref{lem:count}
yields
\[
  \Sigma_T=\Sigma_{T'}\Sigma_{T''}+\Sigma_{T''}\Sigma_{T'} .
\]
If $T'$ or $T''$ is not good, then $T$ is not good and $\Sigma_T=0$ by
induction.  Otherwise $\Sigma_{T'}=2^{|A'|-1}B_{T'}$ and
$\Sigma_{T''}=2^{|A''|-1}B_{T''}$ by induction, and
$B_{T''}B_{T'}=(-1)^{\form{b_{A'}}{b_{A''}}}B_{T'}B_{T''}$ by
Lemma~\ref{lem:calculus}(ii).  Hence
\[
  \Sigma_T=2^{|A|-2}\bigl(1+(-1)^{\form{b_{A'}}{b_{A''}}}\bigr)B_{T'}B_{T''}.
\]
If $\form{b_{A'}}{b_{A''}}=1$, then $T$ is not good and $\Sigma_T=0$.  If
$\form{b_{A'}}{b_{A''}}=0$, then $T$ is good, the matrices $B_{T'}$ and
$B_{T''}$ commute, every ordering inducing $T$ has the product
$B_{T'}B_{T''}$, and $\Sigma_T=2^{|A|-1}B_{T'}B_{T''}$.
\end{proof}

For a general label set we use one tree for each power of two in the binary
expansion of its size.  Let $S\subseteq[M]$ have $k\ge 1$ elements, and let
$k=2^{\alpha_1}+\dots+2^{\alpha_r}$ with $\alpha_1>\dots>\alpha_r\ge 0$ be
the binary expansion of $k$; thus $r$ is the binary digit sum $s_2(k)$ of $k$.  A \emph{forest on $S$}
is a sequence $(T_1,\dots,T_r)$ of trees on pairwise disjoint subsets
$S_1,\dots,S_r$ of $S$ with $|S_q|=2^{\alpha_q}$; these subsets partition $S$.  An
ordering $w$ of $S$ \emph{induces} the forest $(T(w_1),\dots,T(w_r))$, where
$w=w_1w_2\cdots w_r$ is the decomposition of $w$ into consecutive segments of
lengths $2^{\alpha_1},\dots,2^{\alpha_r}$.  A forest is \emph{good} if all its trees are
good, and for a good forest $\Phi=(T_1,\dots,T_r)$ we put
$B_\Phi\coloneqq B_{T_1}\cdots B_{T_r}$.  For $k=M$ we have $r=1$, and a
forest on $[M]$ is a tree on $[M]$.

\begin{corollary}\label{cor:forestsum}
Let $S\subseteq[M]$ have $k\ge 1$ elements, and let $r\coloneqq s_2(k)$.  Then
\[
  \sum_{w}B_w=2^{k-r}\sum_{\Phi}B_\Phi ,
\]
where $w$ runs over the orderings of $S$ and $\Phi$ over the good forests on
$S$.  Each $B_\Phi$ is a signed translation of class $b_S$.
\end{corollary}

\begin{proof}
Every ordering of $S$ induces exactly one forest.  The orderings inducing a
forest $(T_1,\dots,T_r)$ are the concatenations $w_1\cdots w_r$ in which $w_q$
induces $T_q$ for every $q$, so the sum of $B_w$ over these orderings is
$\Sigma_{T_1}\cdots\Sigma_{T_r}$.  By
Lemma~\ref{lem:treesum} this product is $0$ if the forest is not good, and it
equals $2^{(2^{\alpha_1}-1)+\dots+(2^{\alpha_r}-1)}B_\Phi=2^{k-r}B_\Phi$ if
it is good.
\end{proof}

\subsection{Valuations of the logarithmic coefficients}\label{subsec:valuations}

We recall Legendre's formula $\nut(m!)=m-s_2(m)$ for every integer $m\ge 0$.

\begin{proof}[Proof of Proposition~\ref{prop:proper}]
Let $k\coloneqq|S|$ and $r\coloneqq s_2(k)$, so that $1\le k<M$; we have to
show that $\nut(L_S)\ge k$.  By
Lemma~\ref{lem:trace} and Corollary~\ref{cor:forestsum},
\[
  L_S=\frac{(-1)^{k-1}}{2k}\,2^{k-r}\sum_{\Phi}\tr(B_\Phi),
\]
and each trace is $0$, $N$ or $-N$ by Lemma~\ref{lem:calculus}(iii).  Since
$N=2^n$, we obtain
\[
  \nut(L_S)\ge(k-r)+n-1-\nut(k).
\]
By Legendre's formula,
\[
  k-r-\nut(k)=\nut(k!)-\nut(k)=\nut((k-1)!)=k-1-s_2(k-1),
\]
so that $\nut(L_S)\ge k+\bigl(n-2-s_2(k-1)\bigr)$.  Finally
$0\le k-1\le 2^{n-1}-2$.  An integer between $0$ and $2^{n-1}-1$ has at most
$n-1$ binary digits, and its binary digit sum equals $n-1$ only if the integer
is $2^{n-1}-1$.  Hence $s_2(k-1)\le n-2$ and $\nut(L_S)\ge k$.
\end{proof}

For $k=M$ the same computation gives only $\nut(L_{[M]})\ge M-1$, because
$s_2(M-1)=n-1$.  The next proposition expresses $L_{[M]}$ through the good trees on $[M]$.

\begin{proposition}\label{prop:fullformula}
Assume that $d_1+\dots+d_M=0$.  Then $b_{[M]}=(0,0)$.  Consequently, for every
good tree $T$ on $[M]$ there is a sign $\varepsilon_T\in\{1,-1\}$ with
$B_T=\varepsilon_T I$, and
\[
  L_{[M]}=(-1)^{M-1}\,2^{M-1}\sum_{T}\varepsilon_T ,
\]
where $T$ runs over the good trees on $[M]$.  In particular,
$\nut(L_{[M]})=M-1$ if and only if the number of good trees on $[M]$ is odd.
\end{proposition}

\begin{proof}
By~\eqref{eq:classes},
\[
  b_{[M]}=\Bigl(Me_n+\sum_{i=1}^{M}d_i,\ Me_n^*+\sum_{i=1}^{M}\eta_i\Bigr)=(0,0),
\]
because $M$ is even, $d_1+\dots+d_M=0$ by assumption, and
$\eta_1+\dots+\eta_M=0$ by~\eqref{eq:balance}.  For a good tree $T$ on $[M]$
the matrix $B_T$ is therefore a signed translation of class $(0,0)$, that is,
$B_T=\varepsilon_TI$ with $\varepsilon_T\in\{1,-1\}$, and $\tr(B_T)=\varepsilon_TN$.
Lemma~\ref{lem:trace} and Corollary~\ref{cor:forestsum} with $k=M$ and $r=1$
give
\[
  L_{[M]}=\frac{(-1)^{M-1}}{2M}\,2^{M-1}N\sum_{T}\varepsilon_T
  =(-1)^{M-1}\,2^{M-1}\sum_{T}\varepsilon_T ,
\]
since $N=2M$.  The sum $\sum_T\varepsilon_T$ has the same parity as the number
of good trees.
\end{proof}

\begin{example}\label{ex:trees}
In Example~\ref{ex:running} we have $J_0=J(001,e_3^*)$; put
$\lambda\coloneqq e_2^*+e_3^*$, so that $\lambda(x)=x_1+x_2+x_3$.
Lemma~\ref{lem:calculus}(i) gives
\[
  B_1=-W(101,\lambda),\qquad B_2=-W(011,\lambda),\qquad B_3=I,\qquad
  B_4=W(110,0),
\]
so the classes~\eqref{eq:classes} are $b_1=(101,\lambda)$, $b_2=(011,\lambda)$,
$b_3=(000,0)$ and $b_4=(110,0)$.  There are three trees on $\{1,2,3,4\}$,
namely $\{\{1,2\},\{3,4\}\}$, $\{\{1,3\},\{2,4\}\}$ and $\{\{1,4\},\{2,3\}\}$,
where each tree is written as the pair of its cherries.  All three are good;
for instance $\form{b_1}{b_2}=\lambda(101+011)=\lambda(110)=0$.  For the first
tree, $B_1B_2=W(110,0)$ and $B_3B_4=W(110,0)$ by Lemma~\ref{lem:calculus}(i),
so $B_T=W(110,0)^2=I$ and $\varepsilon_T=1$; the other two trees also have
$\varepsilon_T=1$.  Hence $L_{[4]}=-8\cdot 3=-24$ and $\nut(L_{[4]})=3$.

By Lemmas~\ref{lem:calculus}(iii) and~\ref{lem:trace}, $L_S=0$ unless
$b_S=(0,0)$; among the nonempty proper subsets of $[4]$ this happens only for
$\{3\}$ and $\{1,2,4\}$, and $L_{\{3\}}=4$, $L_{\{1,2,4\}}=8$.  Hence
Lemma~\ref{lem:partitions} reads
$[t_{[4]}]F=L_{[4]}+L_{\{3\}}L_{\{1,2,4\}}=-24+32=8$, in agreement with
$Z=\pm 8$.
\end{example}

\section{Parity of the number of good trees}\label{sec:parity}

By Proposition~\ref{prop:fullformula}, Proposition~\ref{prop:full} follows
from the next statement, which holds for every list of nonzero vectors with a
fixed even flag and adapted basis.

\begin{proposition}\label{prop:parity}
The number of good trees on $[M]$ is odd.
\end{proposition}

Since the flag is even, the labels of each height can be grouped into pairs.
We show that the number of good trees keeps its parity when one of these pairs
is required to be a cherry, and we impose the requirements one pair at a time,
in the order of nondecreasing height.  With this order, when a pair of height $j$
is required to be a cherry, every label of smaller height already lies in one
of the cherries required earlier.
Every tree in which all the pairs are cherries is good, and
Lemma~\ref{lem:count} gives the number of these trees.

We recall the definitions used in this section.  The label $i$ has height
$h(d_i)$ and functional $\eta_i=e^*_{h(d_i)}$, and its class is
$b_i=(e_n+d_i,\ e_n^*+\eta_i)$ by~\eqref{eq:classes}; for $S\subseteq[M]$ we
have $b_S=\sum_{i\in S}b_i$.  A tree $T$ on $[M]$ is good if
$\form{b_{A'}}{b_{A''}}=0$ for every node of $T$ of size at least two, with
children $A'$ and $A''$, where
$\form{(a,\beta)}{(c,\gamma)}=\beta(c)+\gamma(a)$.  In particular,
\begin{equation}\label{eq:puredisp}
  \form{(\delta,0)}{(a,\beta)}=\beta(\delta)
  \qquad\text{for all }\delta,a\in V\text{ and }\beta\in V^* .
\end{equation}
We begin with the properties of two labels of equal height.

\begin{lemma}\label{lem:pair}
Let $u\ne v$ be labels of the same height $j$, and put $\delta\coloneqq d_u+d_v$.
\begin{enumerate}[label=\textup{(\roman*)},leftmargin=*]
\item $\delta\in H_{j-1}$ and $b_u+b_v=(\delta,0)$.
\item $\form{(\delta,0)}{b_i}=0$ for every label $i$ of height at least $j$.
\item $\form{(\delta,0)}{b_i+b_{i'}}=0$ for any two distinct labels $i,i'$ of
equal height.
\item $\form{b_u}{b_v}=0$.
\end{enumerate}
\end{lemma}

\begin{proof}
(i) The first statement is~\eqref{eq:sameheight}, applied to the vectors $d_u$
and $d_v$.  Since $\eta_u=\eta_v=e_j^*$, the definition~\eqref{eq:classes}
gives $b_u+b_v=(d_u+d_v,0)$.
(ii) Let $j'\ge j$ be the height of $i$.  The second component of $b_i$ is
$e_n^*+e_{j'}^*$.  Both $e_n^*$ and $e_{j'}^*$ vanish on $H_{j-1}$
by~\eqref{eq:vanish}, and $\delta\in H_{j-1}$, so the claim follows
from~\eqref{eq:puredisp}.
(iii) By~(i), applied to $i$ and $i'$, the second component of $b_i+b_{i'}$ is
zero; now use~\eqref{eq:puredisp}.
(iv) By~(i), $\form{b_u}{b_v}=\form{b_u}{b_u}+\form{b_u}{(\delta,0)}$.  The
first term vanishes because $\form{\kappa}{\kappa}=0$ for every $\kappa$, and
the second vanishes by~(ii), applied to the label $u$.
\end{proof}

For a tree $T$ on $[M]$ and two labels $u\ne v$, let $T^{uv}$ denote the tree
obtained from $T$ by exchanging $u$ and $v$: its nodes are the sets $\tau(A)$,
where $A$ is a node of $T$ and $\tau$ is the transposition of $u$ and $v$.

\begin{lemma}\label{lem:transposition}
Let $T$ be a tree on $[M]$ and let $u\ne v$ be labels.  Then $T^{uv}=T$ if and
only if $\{u,v\}$ is a cherry of $T$.
\end{lemma}

\begin{proof}
If $\{u,v\}$ is a cherry of $T$, then every node of $T$ other than $\{u\}$ and
$\{v\}$ contains both $u$ and $v$ or neither of them, and $\tau$ exchanges
$\{u\}$ and $\{v\}$; hence $T^{uv}$ and $T$ have the same nodes.  Conversely,
suppose that $\{u,v\}$ is not a cherry of $T$, and let $\{u,i\}$ be the cherry
of $T$ containing $u$.  Then $i\ne v$, and $\{v,i\}$ is a node of $T^{uv}$.
It is not a node of $T$, because the only cherry of $T$ containing $i$ is
$\{u,i\}$.  Hence $T^{uv}\ne T$.
\end{proof}

We now fix a partition of $[M]$ into pairs $P_1,\dots,P_{M/2}$, each consisting
of two labels of the same height, and number the pairs so that their heights
are nondecreasing; such a partition exists because the flag is even.  For
$0\le\ell\le M/2$ let $\cG_\ell$ be the set of good trees on $[M]$ in which
each of $P_1,\dots,P_\ell$ is a cherry.  Thus $\cG_0$ is the set of all good
trees on $[M]$.

\begin{lemma}\label{lem:step}
Let $1\le\ell\le M/2$ and $P_\ell=\{u,v\}$.  The map $T\mapsto T^{uv}$ is an
involution of $\cG_{\ell-1}$, and its fixed points are the elements of
$\cG_\ell$.  Consequently $|\cG_{\ell-1}|\equiv|\cG_\ell|\pmod 2$.
\end{lemma}

\begin{proof}
Let $\tau$ be the transposition of $u$ and $v$. Let $T\in\cG_{\ell-1}$,
let $j$ be the height of $u$ and $v$, and put
$\delta\coloneqq d_u+d_v$.  The pairs $P_1,\dots,P_{\ell-1}$ do not contain $u$
or $v$, so they are cherries of $T^{uv}$.  We show that $T^{uv}$ is good.

For every node $C$ of $T$, the elements $b_C$ and $b_{\tau(C)}$ of
$V\times V^*$ differ by $0$ or $(\delta,0)$, because $b_u+b_v=(\delta,0)$ by
Lemma~\ref{lem:pair}(i).  We therefore first show that
\begin{equation}\label{eq:claim}
  \form{(\delta,0)}{b_C}=0
\end{equation}
for every node $C$ of $T$ that does not consist of a single label of height
less than $j$.  (For a single label $i$ of height less than $j$ the left-hand
side of~\eqref{eq:claim} equals $e^*_{h(d_i)}(\delta)$, which need not vanish;
the next paragraph shows that~\eqref{eq:claim} is not needed for such a
node.)  If $C$ consists of a single label of height at least $j$,
then~\eqref{eq:claim} is Lemma~\ref{lem:pair}(ii).  Let $|C|\ge 2$.  Every label of height less than
$j$ belongs to one of the pairs $P_1,\dots,P_{\ell-1}$, because the heights of
the pairs are nondecreasing and $P_\ell$ has height $j$.  These pairs are
cherries of $T$, and a node of size at least two that contains a label
contains the cherry of that label, by the chain property stated after
Definition~\ref{def:tree}.  Hence the labels of height less than $j$
in $C$ form a union of pairs $P_s$ with $s<\ell$.  Each such pair contributes
zero to $\form{(\delta,0)}{b_C}$ by Lemma~\ref{lem:pair}(iii), and each label
of height at least $j$ contributes zero by Lemma~\ref{lem:pair}(ii).  This
proves~\eqref{eq:claim}.

Now let $A$ be a node of $T$ of size at least two, with children $A'$ and
$A''$.  The children of the node $\tau(A)$ of $T^{uv}$ are $\tau(A')$ and
$\tau(A'')$.  By Lemma~\ref{lem:pair}(i) we have
$b_{\tau(A')}=b_{A'}+c'(\delta,0)$, where $c'$ is $1$ if $A'$ contains exactly
one of $u,v$ and $0$ otherwise, and similarly
$b_{\tau(A'')}=b_{A''}+c''(\delta,0)$.  Since
$\form{(\delta,0)}{(\delta,0)}=0$,
\[
  \form{b_{\tau(A')}}{b_{\tau(A'')}}
  =\form{b_{A'}}{b_{A''}}
  +c'\,\form{(\delta,0)}{b_{A''}}
  +c''\,\form{(\delta,0)}{b_{A'}} .
\]
Suppose that $c'=1$.  If $A''$ were a single label $i$ of height less than
$j$, then $A'$ would be a single label as well, namely $u$ or $v$, and $A$
would be a cherry of $T$ containing $i$ and one of $u,v$.  This is impossible,
because the cherry of $T$ containing $i$ is one of the pairs
$P_1,\dots,P_{\ell-1}$.  Hence~\eqref{eq:claim} applies to $C=A''$ and
$\form{(\delta,0)}{b_{A''}}=0$.  In the same way $c''=1$ implies
$\form{(\delta,0)}{b_{A'}}=0$.  Therefore
$\form{b_{\tau(A')}}{b_{\tau(A'')}}=\form{b_{A'}}{b_{A''}}=0$ for every node,
and $T^{uv}$ is good.

Thus $T\mapsto T^{uv}$ maps $\cG_{\ell-1}$ to itself, and it is an involution.
By Lemma~\ref{lem:transposition} its fixed points are the trees of
$\cG_{\ell-1}$ in which $P_\ell$ is a cherry, that is, the elements of
$\cG_\ell$.  The remaining elements of $\cG_{\ell-1}$ are exchanged in pairs.
\end{proof}

In Example~\ref{ex:trees} the pairs are $P_1=\{1,2\}$ and $P_2=\{3,4\}$.  The
exchange of the labels $1$ and $2$ fixes the tree $\{\{1,2\},\{3,4\}\}$ and
exchanges the other two trees, in accordance with Lemma~\ref{lem:step}.

\begin{lemma}\label{lem:last}
Every tree on $[M]$ in which each of $P_1,\dots,P_{M/2}$ is a cherry is good.
Consequently $|\cG_{M/2}|=(M/2)!/2^{M/2-1}$, and this number is odd.
\end{lemma}

\begin{proof}
Let $T$ be such a tree.  Its cherries are exactly $P_1,\dots,P_{M/2}$, since
the cherries of a tree on $[M]$ partition $[M]$.  Let $A$ be a node of $T$ of
size at least two with children $A'$ and $A''$.  If $|A|=2$, then
$A=P_s=\{u,v\}$ for some $s$, and
$\form{b_u}{b_v}=0$ by Lemma~\ref{lem:pair}(iv).  If $|A|\ge 4$, then $A'$ and
$A''$ are unions of pairs $P_s$, so $b_{A'}=(\delta',0)$ and
$b_{A''}=(\delta'',0)$ for some $\delta',\delta''\in V$ by
Lemma~\ref{lem:pair}(i), and $\form{b_{A'}}{b_{A''}}=0$
by~\eqref{eq:puredisp}.  Hence $T$ is good.

Therefore $\cG_{M/2}$ is the set of all trees on $[M]$ whose cherries are
$P_1,\dots,P_{M/2}$.  Replacing each cherry by a single element gives a
bijection between this set and the set of trees on the set
$\{P_1,\dots,P_{M/2}\}$, which has $m\coloneqq M/2=2^{n-2}$ elements.  By
Lemma~\ref{lem:count} their number is $m!/2^{m-1}$, and it is odd because
$\nut(m!)=m-1$ by Legendre's formula.
\end{proof}

\begin{proof}[Proof of Proposition~\ref{prop:parity}]
By Lemma~\ref{lem:step},
$|\cG_0|\equiv|\cG_1|\equiv\dots\equiv|\cG_{M/2}|\pmod 2$, and
$|\cG_{M/2}|$ is odd by Lemma~\ref{lem:last}.
\end{proof}

\begin{proof}[Proof of Proposition~\ref{prop:full}]
Let $d_1+\dots+d_M=0$.  The number of good trees on $[M]$ is odd by
Proposition~\ref{prop:parity}, so $\nut(L_{[M]})=M-1$ by
Proposition~\ref{prop:fullformula}.
\end{proof}

For a zero-sum list, the one-block term in Lemma~\ref{lem:partitions}
is therefore the unique term of valuation $M-1$. By~\eqref{eq:ZF},
$\nut(Z)=M-1$. In particular $Z\ne0$, and Lemma~\ref{lem:certificate}
gives a perfect matching with the prescribed differences.

\section{Concluding remarks}\label{sec:conclusion}

We have proved the conjecture of Balister, Gy\H{o}ri and Schelp: for $n\ge 2$,
a list of $M=2^{n-1}$ nonzero vectors of $\F^n$ is the list of differences of a
perfect matching of $\F^n$ exactly when its sum is zero.  The theorem applies
to every such list, whatever the multiplicities of its vectors and their
position relative to the hyperplanes of $\F^n$.  The proof gives the stronger
Theorem~\ref{thm:valuation}: for a list with sum zero and the signs defined by
any total order of $V$ and any basis of $V$ adapted to an even flag for the list, the
signed count of realizations is an odd multiple of $2^{M-1}$.

The argument uses linear algebra over $\F$, the expansion~\eqref{eq:pfaffian}
of the Pfaffian together with the identities~\eqref{eq:pf-facts}, and formal
logarithms of the normalized Pfaffian. In the proof of sufficiency,
the zero-sum condition makes the class $b_{[M]}$ vanish in
Proposition~\ref{prop:fullformula}.
The even flag, which exists for every list of $M$ nonzero vectors
(Lemma~\ref{lem:flag}), supplies the relation $\eta_1+\dots+\eta_M=0$ and the
pairs of labels of equal height on which the parity argument of
Section~\ref{sec:parity} rests.

The proof of Theorem~\ref{thm:main} is existential: it shows that the
realizations cannot cancel completely in the signed count.  An efficient
construction of a realization from the list $(d_1,\dots,d_M)$ would be of
interest.


\begin{thebibliography}{99}

\bibitem{BGS}
P.~N. Balister, E.~Gy\H{o}ri and R.~H. Schelp,
Coloring vertices and edges of a graph by nonempty subsets of a set,
\emph{European J. Combin.} \textbf{32} (2011), no.~4, 533--537.

\bibitem{Godsil}
C.~D. Godsil,
\emph{Algebraic Combinatorics},
Chapman \& Hall, New York, 1993.

\bibitem{GK}
L.~Golowich and C.~Kim,
New classes of set-sequential trees,
\emph{Discrete Math.} \textbf{343} (2020), no.~3, 111741,
\href{https://doi.org/10.1016/j.disc.2019.111741}{doi:10.1016/j.disc.2019.111741}.

\bibitem{Hall}
M.~Hall, Jr.,
A combinatorial problem on abelian groups,
\emph{Proc. Amer. Math. Soc.} \textbf{3} (1952), 584--587.

\bibitem{HKRS}
H.~D.~L. Hollmann, K.~Khathuria, A.-E. Riet and V.~Skachek,
On some batch code properties of the simplex code,
\emph{Des. Codes Cryptogr.} \textbf{91} (2023), no.~5, 1595--1605,
\href{https://doi.org/10.1007/s10623-022-01173-6}{doi:10.1007/s10623-022-01173-6}.

\bibitem{Kovacs}
B.~Kov\'acs,
Finding a perfect matching of $\mathbb{F}_2^n$ with prescribed differences,
\emph{Ars Math. Contemp.} \textbf{26} (2026), no.~1, P1.05,
\href{https://doi.org/10.26493/1855-3974.3265.91b}{doi:10.26493/1855-3974.3265.91b}.

\bibitem{FourEight}
A.~Kreindel and A.~L. Zabokritskiy (Yohananov),
Prescribed-difference matchings with four and eight holes: Fourier filters for
compatible boundaries,
\href{https://openalex.org/works/W7204825102}{preprint} 2608.27547v1, 2026.

\bibitem{YE}
L.~Yohananov and I.~B. Essayag,
Optimal functional $2^{s-1}$-batch codes: exploring new sufficient conditions,
in \emph{2025 IEEE International Symposium on Information Theory (ISIT)},
2025, 1--6, \href{https://doi.org/10.1109/ISIT63088.2025.11195523}{doi:10.1109/ISIT63088.2025.11195523}.

\bibitem{YY}
L.~Yohananov and E.~Yaakobi,
Almost optimal construction of functional batch codes using extended Simplex
codes,
\emph{IEEE Trans. Inform. Theory} \textbf{68} (2022), no.~10, 6434--6451,
\href{https://doi.org/10.1109/TIT.2022.3162158}{doi:10.1109/TIT.2022.3162158}.

\bibitem{TwoHole}
A.~L. Zabokritskiy (Yohananov),
Perfect matchings with prescribed differences beyond Hall: the two-hole
problem,
\href{https://openalex.org/works/W7167968659}{preprint} 2607.08630v1, 2026.

\bibitem{ZEY}
Y.~Zhang, T.~Etzion and E.~Yaakobi,
Bounds on the length of functional PIR and batch codes,
\emph{IEEE Trans. Inform. Theory} \textbf{66} (2020), no.~8, 4917--4934,
\href{https://doi.org/10.1109/TIT.2020.2977631}{doi:10.1109/TIT.2020.2977631}.

\end{thebibliography}
\end{document}